\documentclass[12pt,reqno]{article}

\usepackage[usenames]{color}
\usepackage{amssymb}
\usepackage{graphicx}
\usepackage{amscd}

\usepackage[colorlinks=true,
linkcolor=webgreen,
filecolor=webbrown,
citecolor=webgreen]{hyperref}

\definecolor{webgreen}{rgb}{0,.5,0}
\definecolor{webbrown}{rgb}{.6,0,0}

\usepackage{color}
\usepackage{fullpage}
\usepackage{float}

\usepackage{psfrag}
\usepackage{graphics,amsmath,amssymb}
\usepackage{amsthm}
\usepackage{amsfonts}
\usepackage{latexsym}
\usepackage{epsf}

\usepackage{enumitem}
\usepackage[mathscr]{euscript}

\DeclareMathOperator{\Hom}{Hom}
\DeclareMathOperator{\bij}{bij}
\DeclareMathOperator{\inj}{inj}
\DeclareMathOperator{\sur}{sur}

\begin{document}

\theoremstyle{plain}
\newtheorem{theorem}{Theorem}
\newtheorem{corollary}[theorem]{Corollary}
\newtheorem{lemma}[theorem]{Lemma}
\newtheorem{proposition}[theorem]{Proposition}

\theoremstyle{definition}
\newtheorem{definition}[theorem]{Definition}
\newtheorem{example}[theorem]{Example}
\newtheorem{conjecture}[theorem]{Conjecture}

\theoremstyle{remark}
\newtheorem{remark}[theorem]{Remark}

\begin{center}
\vskip 1cm{\LARGE\bf Partial Profiles of Quasi-Complete Graphs
}
\vskip 1cm
\large
Pedro Lopes\\
Center for Analysis, Geometry, and Dynamical Systems\\
Department of Mathematics\\
Instituto Superior T\'ecnico\\
University of Lisbon\\
Av. Rovisco Pais\\
1049-001 Lisbon\\
Portugal\\
\href{mailto:pelopes@math.tecnico.ulisboa.pt}{\tt pelopes@math.tecnico.ulisboa.pt}
\end{center}

\vskip .2 in

\begin{abstract}
We enumerate graph homomorphisms to quasi-complete graphs, i.e.,
graphs obtained from complete graphs by removing one edge. The source
graphs are complete graphs, quasi-complete graphs, cycles, paths, wheels and broken wheels. These enumerations give rise to sequences of integers with two indices; one of the indices is the number of vertices of the source graph, and the other index is the number of vertices of the target graph.
\end{abstract}

\section{Introduction}
\label{sec:introduction}

We believe that the enumeration of graph homomorphisms \cite{hell} to a given graph is an interesting problem in itself to which not much attention has been devoted -- except perhaps for the chromatic polynomials of a few families of graphs \cite{dkt}. In this article we take up the task of calculating the numbers of homomorphisms to the so-called quasi-complete graphs, i.e., graphs obtained from complete graphs by removing one edge. We take for source graphs, complete graphs, quasi-complete graphs, paths, cycles, wheels and broken wheels. Specifically, we provide formulas for the numbers of these homomorphisms in terms of the numbers of vertices of source and target graphs.

The availability of such formulas may aid in the characterization and/or the identification of graphs \cite{Dvorak}, leaning on Lov\'asz' theorem \cite{ll}.  Furthermore, it may help supporting or providing counter-examples to the reconstruction conjecture \cite{borgs}. In fact, `` \ldots the reconstruction conjecture is equivalent to the assertion that it is enough to know all numbers $\hom (F, G)$ with $|V(F)|<|V(G)|$ in order to recover the isomorphism type of $G$.'' \cite[1st paragraph,  p.\ 317]{borgs}. In this context, having a catalog of numbers $\hom(F, G)$, for a significant number of pairs of graphs $(F, G)$, could lead to finding a pair of non-isomorphic graphs, $G$ and $G'$, for which $\hom(F, G) = \hom(F, G')$ for all $F$ with $|V(F)|<|V(G)|$. In this case there would be a counter-example to the reconstruction conjecture. Otherwise, it would be supporting evidence for the conjecture.

We remark that other aspects of graph homomorphisms are currently the object of  research. These include, but do not exhaust, topics like extremal issues in graph homomorphisms \cite{LohPikhurkoSudakov, CsikvariLin1, CsikvariLin2, Galvin} or the complexity of enumeration of graph homomorphisms \cite{DyerGreenhill}. We emphasize that the present contribution is to the exact  enumeration of these homomorphisms.

We now recall the definitions of graph,  of graph homomorphism, and of other objects used here in order to develop the notation and the terminology. A graph is a finite set of points, called \emph{vertices}, along with a set of unordered pairs of distinct vertices, called \emph{edges}. If there is an edge joining vertices $i$ and $j$, it is denoted $ij$. In this case, we say that vertices $i$ and $j$ are \emph{adjacent}. A graph is denoted by capital letters $G$ or $H$. The set of vertices of graph $G$ is denoted $V(G)$ and the set of edges of $G$ is denoted $E(G)$.  A \emph{homomorphism}, $f$, from a graph $G$ to a graph $H$, is a map from $V(G)$ to $V(H)$ which preserves adjacency, i.e., whenever $ij$ is an edge of $G$, $f(i)f(j)$ is an edge of $H$ \cite{hell}. The set of homomorphisms from $G$ to $H$ is denoted $\Hom(G, H)$ and its cardinality $\hom(G, H)$. Also, $\bij (G, H)$, $\sur (G, H)$ and $\inj (G, H)$ denote the cardinality of the bijective homomorphisms, the surjective homomorphisms and the injective homomorphisms of $\Hom (G, H)$.

Let $(H\sb{1}, H\sb{2}, H\sb{3}, \ldots )$ be a fixed sequence of all non-isomorphic graphs. Given a graph $G$, the \emph{profile} of $G$ (or \emph{Lov\'asz vector} of $G$) is the sequence of non-negative integers $(n\sb{1}, n\sb{2}, n\sb{3}, \ldots )$, where $n\sb{i}=\hom(H\sb{i}, G)$. Lov\'asz \cite{ll} introduced the profile of a graph in order to prove that graphs with the same profile are
isomorphic modulo twins. We use the term \emph{partial profile} of a given graph since we only calculate a subsequence (albeit infinite) of the $n\sb{i}$'s.

For each graph $G$, the chromatic polynomial of $G$, $\chi\sb{G} (m)$, is the number of homomorphisms from $G$ to the complete graph, $K\sb{m}$, for each $m$ greater than or equal to the chromatic number of $G$, $\chi\sb{G}(m)=\hom(G, K\sb{m})$ \cite{dkt}. As such, all the known chromatic polynomials evaluated at a specific positive integer $m$, provide a partial profile of the complete graph $K\sb{m}$.

The rest of this article is organized as follows. In Section \ref{sec:count} we prove the formulas we obtained for the numbers of graph homomorphisms. The target graphs are always quasi-complete graphs; except for Subsection \ref{subsect:aux}, each subsection in Section \ref{sec:count} is devoted to a particular kind for source graph. In Subsection \ref{subsect:aux} we provide auxiliary formulas to help in the calculations of the subsequent enumerations. In Section \ref{finrems} we formulate a few questions for future work. In Section \ref{ack} we acknowledge our financial sponsors, hosts, and colleagues.

\section{Counting homomorphisms to quasi-complete graphs}\label{sec:count}

For any integer $m\geq 3$, we let $K\sb{m}$ denote the \emph{complete graph on $m$
vertices}, i.e., the graph on $m$ vertices such that any two vertices are adjacent. For any integer $m\geq 3$, we define the \emph{quasi-complete graph} on $m$ vertices to be the graph obtained from $K\sb{m}$ by removing one edge. We denote it $K\sb{m}\sp{1}$. Furthermore, the edge that has been removed in order to obtain it from $K\sb{m}$ is referred to as the \emph{exceptional edge}. The vertices which make up this edge are denoted $A$ and $B$ throughout the article.

In this section, we obtain formulas for calculating the numbers of homomorphisms into quasi-complete graphs from certain graphs in terms of the numbers of vertices. Except for Subsection \ref{subsect:aux}, the title of each subsection identifies which class of graphs is being used as source graph.

\subsection{The source graphs are complete graphs}\label{subsect:complete}

\begin{proposition}For integers $n\geq 3$, $m\geq 3$,
\begin{displaymath}
\hom(K\sb{n}, K\sb{m}\sp{1}) = \begin{cases}
0, &  \text{if  $n\geq m$;}\\
(m-1)!\cdot 2, & \text{if  $n=m-1$;}\\
n!{{m-2}\choose {n-1}}\cdot 2 + {{m-2}\choose {n}} n!, & \text{if  $n\leq m-2$.}
\end{cases}
\end{displaymath}
In particular,
\[
\hom(K\sb{n}, K\sb{m}\sp{1}) = \inj(K\sb{n}, K\sb{m}\sp{1}) , \qquad
\sur(K\sb{n}, K\sb{m}\sp{1}) = \bij(K\sb{n}, K\sb{m}\sp{1}) = 0  ,
\]
for all $n \leq m-1$.
\end{proposition}
\begin{proof} Homomorphisms from complete graphs to any given
graph are one-to-one since any two vertices in a complete graph
are adjacent. In this way, if $n > m$, then, there are no
homomorphisms from $K\sb{n}$ to $K\sb{m}\sp{1}$.

If $n=m$, then any
one-to-one map from $K\sb{n}$ to $K\sb{m}\sp{1}$ has to
send distinct vertices of $K\sb{n}$ to the vertices of the exceptional edge of
$K\sb{m}\sp{1}$, thereby mapping adjacent vertices to non-adjacent
vertices. Again, there are no homomorphisms in this case.

If $n=m-1$, then the image of a homomorphism cannot contain both vertices of the exceptional edge as pointed out above. On the other hand, it has to contain at least one of these two vertices, for otherwise there would not be enough vertices in the image. The problem is then reduced to distribute the $n$ vertices of $K\sb{n}$, in a one-to-one fashion, over the $m-1(=n)$ vertices of $K\sb{m}\sp{1}$ which contain exactly one of the vertices of the exceptional edge, or over those that contain the other vertex of the exceptional edge. This is
\[
n!\cdot 2 = (m-1)!\cdot 2 .
\]

Now, for $n<m-1$. As in the preceding case, the image of a homomorphism may contain exactly one vertex of the exceptional edge. In this case, after selecting one of the vertices from the exceptional edge, the $n$ vertices of $K\sb{n}$ are distributed over a subset of $n$ vertices from the set of $m-1$ vertices of $K\sb{m}\sp{1}$ which contain the selected vertex of the exceptional edge (but not the other one). Then, one of the vertices of $K\sb{n}$ has to be mapped to this vertex ($n$ possibilities). The remaining $n-1$ vertices of $K\sb{n}$ have to be mapped, in a one-to-one fashion, to the vertices of $K\sb{m}\sp{1}$ but the vertices of the exceptional edge, to ensure adjacency is preserved (${{m-2}\choose{n-1}}(n-1)!$ possibilities). The other vertex of the exceptional edge is then selected and the preceding reasoning is repeated yielding
\[
2n{{m-2}\choose{n-1}}(n-1)!=2 n!{{m-2}\choose{n-1}} .
\]

Since $n<m-1$, it is still possible to construct homomorphisms whose image do not contain either of the vertices of the exceptional edge. In this case, the $n$ vertices of $K\sb{n}$ are mapped in a one-to-one fashion to the $m-2$ vertices of $K\sb{m}\sp{1}$, except the vertices from the exceptional edge. There are
\[
{{m-2}\choose{n}}n!
\]
such possibilities.

Clearly, the numbers $\inj (K\sb{n}, K\sb{m}\sp{1})$ coincide with the numbers $\hom (K\sb{n}, K\sb{m}\sp{1})$ since the source graphs are complete graphs. For the same reason, $\sur (K\sb{n}, K\sb{m}\sp{1}) = \bij (K\sb{n}, K\sb{m}\sp{1}) = 0$.
\end{proof}

\subsection{The source graphs are quasi-complete graphs}\label{subsect:sourceqc}

\begin{proposition}For integers $n\geq 3$, $m\geq 3$,
\begin{displaymath}
\hom(K\sb{n}\sp{1}, K\sb{m}\sp{1}) = \\
\begin{cases}
0, & \text{if $n \geq  m+1$;}\\
2(n-1)! + 2(n-2)!, & \text{if $m = n$;}\\
\displaystyle{\frac{2(m-2)!}{(m-n)!}} + \frac{\bigl( 2(n-1) + 3\bigr) (m-2)!}{(m-n-1)!}   + \frac{(m-2)!}{(m-n-2)!}, & \text{if $m > n$.}
\end{cases}
\end{displaymath}
Also,
\[
\bij(K\sb{n}\sp{1}, K\sb{n}\sp{1}) = \sur(K\sb{n}\sp{1}, K\sb{n}\sp{1}) = \inj(K\sb{n}\sp{1}, K\sb{n}\sp{1}) = (n-2)!2 ,
\]
and
\begin{align*}
&\inj(K\sb{n}\sp{1}, K\sb{m}\sp{1}) \\
&=\displaystyle{\frac{2(m-2)!}{(m-n)!}} + \frac{2(n-1)(m-2)!(m - n-1)}{(m-n)!}   + \frac{2(m-2)!}{(m-n-1)!} +  \frac{(m-2)!}{(m-n-2)!}, \qquad  m > n ,
\end{align*}
and
\[
\bij (K\sb{n}\sp{1}, K\sb{m}\sp{1}) = \sur (K\sb{n}\sp{1}, K\sb{m}\sp{1}) = 0 ,  \qquad \text{ for } m > n .
\]
\end{proposition}
\begin{proof} We regard $K\sb{n}\sp{1}$ as consisting of $K\sb{n-1}$ formed by vertices $1, 2, \ldots , n-1 $ and then adjoining vertex $n$ making it adjacent to all other vertices but $1$. Concerning $K\sb{m}\sp{1}$, we denote the vertices of the exceptional edge $A$ and $B$.

Any homomorphism from $K\sb{n}\sp{1}$ to $K\sb{m}\sp{1}$ has to map the $K\sb{n-1}$ in $K\sb{n}\sp{1}$ to a complete subgraph of $K\sb{m}\sp{1}$ with $n-1$ vertices. It follows that, if $n-1 > m-1$, i.e., if $n > m$, then $\hom (K\sb{n}\sp{1}, K\sb{m}\sp{1}) = 0$.

If $m = n$, then there are two complete subgraphs of $K\sb{m}\sp{1}$ with $n-1$ vertices, one containing vertex $A$, and the other one containing vertex $B$. It follows that the image of each homomorphism, in this $m=n$ case, has to involve    at least one of $A$ and $B$.

The homomorphisms that involve both have to map $1$ to $A$ and $n$ to $B$, or $1$ to $B$ and $n$ to $A$. The other $n-2$ vertices are permuted among the $n-2$ vertices of $K\sb{m}\sp{1}$ other than $A$ and $B$. So the homomorphisms that involve both $A$ and $B$, i.e., the bijective ones, are
\[
\bij (K\sb{n}\sp{1}, K\sb{n}\sp{1}) = \sur (K\sb{n}\sp{1}, K\sb{n}\sp{1}) = \inj (K\sb{n}\sp{1}, K\sb{n}\sp{1}) = (n-2)!2 .
\]

Now, for the homomorphisms that involve exactly one of $A$ or $B$. The $K\sb{n-1}$ of $K\sb{n}\sp{1}$ is mapped to either one of the two complete subgraphs of $K\sb{m}\sp{1}$ with $n-1$ vertices mentioned above, and vertex $n$ is mapped to the same vertex $1$ is mapped to. These are $(n-1)!2$. Then,
\[
\hom (K\sb{n}\sp{1}, K\sb{n}\sp{1}) = (n-1)!2 + (n-2)!2 .
\]

Now, suppose $n < m$. Besides homomorphisms involving both $A$ and $B$, and homomorphisms involving either $A$ or $B$ (exclusively), there are now homomorphisms which do not involve $A$ or $B$.

The ones that involve both $A$ and $B$ are injective. Also, $1$ is mapped to $A$ and $n$ is mapped to $B$, or $1$ is mapped to $B$ and $n$ is mapped to $A$, while the remaining $n-2$ vertices are permuted over the $m-2$ vertices of $K\sb{m}\sp{1}$ other than $A$ and $B$. These are, then,
\[
2{{m-2}\choose{n-2}}(n-2)! .
\]

Now, for the homomorphisms that involve exactly $A$ or $B$. One of the $n$ vertices has to be mapped to $A$ (respect., $B$, and for this reason there will be an overall $2$ factor). Consider the vertices $1, 2, \ldots , n-1$ forming the $K\sb{n-1}$ in $K\sb{n}\sp{1}$. One of these is mapped to $A$ ($n-1$ possibilities). Then the remaining $n-2$ vertices of $K\sb{n-1}$ are distributed over $n-2$ of $m-2$ vertices of $K\sb{m}\sp{1}$ (other than $A$ and $B$). Finally, vertex $n$ can either be mapped to the vertex $1$ has been mapped to, or to one of the remaining $m-2-(n-1) = m-n-1$ vertices of $K\sb{m}\sp{1}$ (these latter homomorphisms are injective). If vertex $n$ is mapped to $A$, then vertex $1$ is not mapped to $A$ (this situation has already been contemplated) and the remaining $n-1$ vertices of $K\sb{n}$, including vertex $1$, are distributed over the remaining $m-2$ vertices of $K\sb{m}\sp{1}$. Thus,
\[
2(n-1){{m-2}\choose{n-2}}(n-2)! \biggl( 1 + (m - n - 1) \biggr) + 2 {{m-2}\choose{n-1}}(n-1)!
\]
is the number of homomorphisms that involve exactly one of $A$ or $B$.
The injective homomorphisms among these are
\[
2(n-1){{m-2}\choose{n-2}}(n-2)!( m - n - 1 ) + 2 {{m-2}\choose{n-1}}(n-1)! .
\]

Finally, those that do not involve $A$ nor $B$ are
\[
{{m-2}\choose{n}}n! + {{m-2}\choose{n-1}}(n-1)! ,
\]
where the left summand counts the injective ones of these homomorphisms, and the right summand counts the homomorphisms which map $1$ and $n$ to the same vertex. Then, for $n<m$,
\begin{multline*}
\inj (K\sb{n}\sp{1}, K\sb{m}\sp{1})  \\
= 2{{m-2}\choose{n-2}}(n-2)! + 2(n-1){{m-2}\choose{n-2}}(n-2)!( m - n - 1 ) + 2 {{m-2}\choose{n-1}}(n-1)! + {{m-2}\choose{n}}n! ,
\end{multline*}
and
\begin{multline*}
\hom (K\sb{n}\sp{1}, K\sb{m}\sp{1}) = 2{{m-2}\choose{n-2}}(n-2)! + 2(n-1){{m-2}\choose{n-2}}(n-2)! (m - n)  \\
 + 2 {{m-2}\choose{n-1}}(n-1)! + {{m-2}\choose{n}}n! + {{m-2}\choose{n-1}}(n-1)! .
\end{multline*}

The proof is concluded.
\end{proof}

\subsection{Auxiliary material}\label{subsect:aux}

This section is intended as auxiliary material for the understanding of the proofs of Propositions \ref{prop:pnkm1} and \ref{prop:cnkm1}.

\begin{figure}[h!]
    \psfrag{0}{\huge $\mathbf 0$}
    \psfrag{1}{\huge $\mathbf 1$}
    \centerline{\scalebox{.50}{\includegraphics{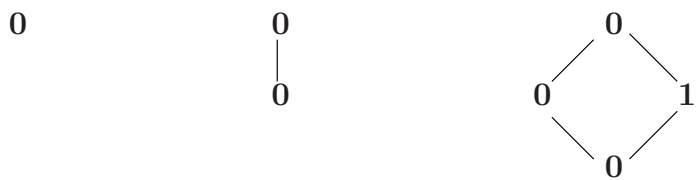}}}
    \caption{Graphs $G\sb{1}$ (left), $G\sb{2}$ (middle), and $G\sb{3}$ (right).}\label{Fi:graphs123}
\end{figure}

\begin{figure}[h!]
    \psfrag{0}{\huge $\mathbf 0$}
    \psfrag{1}{\huge $\mathbf 1$}
    \centerline{\scalebox{.50}{\includegraphics{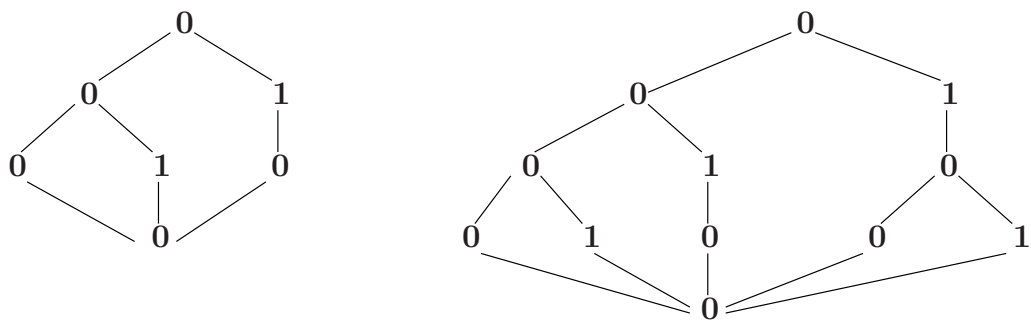}}}
    \caption{Graphs $G\sb{4}$ (left),  and $G\sb{5}$ (right).}\label{Fi:graphs45}
\end{figure}

\begin{definition}\label{def:graphs} For each integer $i \geq 1$, we define the $($directed$)$ graph
 $G\sb{i}$ as follows. The vertices of $G\sb{i}$ are $0$'s
and $1$'s. The vertices are drawn on different
levels and the levels are displayed vertically; higher levels are drawn below lower levels. The graph $G\sb{i}$ has $i$ levels, $1, 2, \ldots , i$. The first level is the \emph{top level}, the last level is the
\emph{bottom level}, the other levels are the \emph{intermediate levels}. The
first level has only one vertex, a $0$. Likewise, the last level has only one vertex,
a $0$. The edges connect vertices on consecutive levels only and
are directed from top to bottom $($and this is why we omit the
arrows on the edges$)$. Graph $G\sb{1}$ is formed by only one vertex, a $0$; graph $G\sb{2}$ is formed by two vertices $($two $0$'s$)$ and an edge connecting them $($see Figure \ref{Fi:graphs123}$)$. For $i>2$, graph $G\sb{i}$ is constructed in the following way. For each $1 < j < i$, each vertex on level $j$ is the endpoint of exactly one edge stemming from level $j-1$.  Up to and including level $i-2$, from each
$0$ stems an edge to a $0$, and an edge to a $1$, whereas from each
$1$ stems only one edge to a
$0$. From each vertex on level $i-1$ stems one edge to the vertex, $0$, on the bottom level.
Graph $G\sb{3}$ is depicted in
Figure \ref{Fi:graphs123}. It
has three levels, the first one with a $0$, the second one with a $0$
and a $1$, and the third one with one $0$. See also Figure
\ref{Fi:graphs45} for graphs $G\sb{4}$ and $G\sb{5}$.
\end{definition}

\begin{definition}\label{def:piqi} Leaning on the definition of the directed graphs $G\sb{i}$, see Definition \ref{def:graphs} above,
we now define polynomials on a variable $m$, the
$p\sb{i}(m)$ and the $q\sb{i}(m)$. By definition,
\[
p\sb{0}\equiv 1 , \qquad q\sb{0}\equiv 1 .
\]
For integer $i\geq 1$, the $p\sb{i}$ and $q\sb{i}$ are to be read off the corresponding graph $G\sb{i}$ as
follows.
Fix integer $i\geq 1$. Each $($maximal$)$ directed path on $G\sb{i}$ $($with respect to the number of vertices$)$, starting at the $0$ on
the top level and ending at the $0$ on the bottom level,
gives rise to a summand of $p\sb{i}$ $($respect., $q\sb{i}$$)$,
formed by the product of $i$ factors among $2$, $(m-2)$, and
$(m-3)$ $($respect., $2$, $(m-2)$, $(m-3)$, and $(m-1)$$)$. Some
of these factors may be missing or may be repeated in a given
summand.

The $0$ on the top level stands for an $(m-2)$ factor.
Each $1$ stands for a $2$ factor. A $0$
stemming from a $1$ stands always for an $(m-2)$ factor. A $0$ on any
intermediate level
and stemming from a $0$ stands for an $(m-3)$ factor; on the bottom level, it stands for an $(m-3)$
for $p\sb{i}$, whereas it stands for an $(m-1)$ for the $q\sb{i}$.
So, for each directed path starting at the $0$ on the top level and ending
at the $0$ on the bottom level, the $0$'s and the $1$'s
stand for the indicated factors, and the edges of the path are to
be considered as an instruction for multiplying the factors connected by the edges.

 The sum of these
products over all maximal directed paths of $G\sb{i}$, starting at the top
$0$ and ending at the bottom $0$, yields the $p\sb{i}$
$($respect., $q\sb{i})$.
\end{definition}

Here are some examples:

\[
p\sb{1}(m) = m-2 , \qquad q\sb{1}(m) = m-2 ;
\]

\[
p\sb{2}(m) = (m-2)(m-3) , \qquad q\sb{2}(m) = (m-2)(m-1) ;
\]

\begin{align*}
p\sb{3}(m) &= (m-2)(m-3)(m-3)+(m-2)2(m-2) , \\
q\sb{3}(m) &= (m-2)(m-3)(m-1)+(m-2)2(m-2) ;
\end{align*}

\begin{multline*}
p\sb{4}(m) = \\
(m-2)(m-3)(m-3)(m-3)+(m-2)(m-3)2(m-2)+(m-2)2(m-2)(m-3) ,
\end{multline*}
\begin{multline*}
q\sb{4}(m) = \\
(m-2)(m-3)(m-3)(m-1)+(m-2)(m-3)2(m-2)+(m-2)2(m-2) ;
\end{multline*}

\begin{multline*}p\sb{5}(m) =
(m-2)(m-3)(m-3)(m-3)(m-3)+(m-2)(m-3)(m-3)2(m-2)+\\
+(m-2)(m-3)2(m-2)(m-3)+
(m-2)2(m-2)(m-3)(m-3)+(m-2)2(m-2)2(m-2) ,
\end{multline*}
\begin{multline*}
q\sb{5}(m) =
(m-2)(m-3)(m-3)(m-3)(m-1)+(m-2)(m-3)(m-3)2(m-2)+\\
+(m-2)(m-3)2(m-2)(m-1)+(m-2)2(m-2)(m-3)(m-1)+(m-2)2(m-2)2(m-2) .
\end{multline*}

\subsection{The source graphs are paths}\label{subsect:path}

Let $n$ be an integer greater than $1$. We remark that a \emph{path on $n$ vertices}, denoted $P\sb{n}$, is a graph with vertices $0, 1, \ldots , n-1$ and edges $01, 12, \ldots , (n-2)(n-1)$.

\bigbreak

\begin{proposition}\label{prop:pnkm1}For integers $n\geq 2$, $m\geq 3$,
\[
\hom(P\sb{n}, K\sb{m}\sp{1}) = m(m-1)\sp{n-1} -
\sum\sb{k=1}\sp{n-1}s\sb{n, m}\sp{k} ,
\]
where
\[
s\sb{n, m}\sp{n-1}  = 2 ,
\]
and, for $1\leq k \leq n-2$,
\begin{align*}
s\sb{n, m}\sp{k} & = \frac{1}{2}\sum\sb{\substack{\{ 0 = i\sb{1}<i\sb{2}< \cdots < i\sb{k} = n-2 \} \\
\text{ s.t. } i\sb{j}-i\sb{j-1}\neq 2}}
 \Biggl[  \prod\sb{\substack{i\sb{j} \, \text{ s.t. }\\
i\sb{j}-i\sb{j-1}>2 \\  2\leq j \leq k}} 2\,
p\sb{i\sb{j}-i\sb{j-1}-2}(m) \Biggr]  \quad  \\
& \qquad + \sum\sb{\substack{\{ 0 = i\sb{1}<i\sb{2}< \cdots < i\sb{k} < n-2 \} \\
\text{ s.t. } i\sb{j}-i\sb{j-1}\neq 2}}
\Biggl[  \prod\sb{\substack{i\sb{j} \, \text{ s.t. }\\
i\sb{j}-i\sb{j-1}>2 \\  2\leq j \leq k}} 2\,
p\sb{i\sb{j}-i\sb{j-1}-2}(m) \Biggr]  \cdot q\sb{n-2-i\sb{k}}(m)  \\
& \qquad + \sum\sb{\substack{\{ 0 < i\sb{1}<i\sb{2}< \cdots < i\sb{k} = n-2 \} \\
\text{ s.t. } i\sb{j}-i\sb{j-1}\neq 2}}q\sb{i\sb{1}}(m)
\cdot \Biggl[  \prod\sb{\substack{i\sb{j} \, \text{ s.t. }\\
i\sb{j}-i\sb{j-1}>2 \\  2\leq j \leq k}} 2\,
p\sb{i\sb{j}-i\sb{j-1}-2}(m) \Biggr]   \\
& \qquad + 2 \sum\sb{\substack{\{ 0 < i\sb{1}<i\sb{2}< \cdots < i\sb{k} < n-2 \} \\
\text{ s.t. } i\sb{j}-i\sb{j-1}\neq 2}}q\sb{i\sb{1}}(m)
\cdot \Biggl[  \prod\sb{\substack{i\sb{j} \, \text{ s.t. }\\
i\sb{j}-i\sb{j-1}>2 \\  2\leq j \leq k}} 2\,
p\sb{i\sb{j}-i\sb{j-1}-2}(m) \Biggr]  \cdot q\sb{n-2-i\sb{k}}(m) ,
\end{align*}
and the $p\sb{i}$'s and $q\sb{i}$'s are the polynomials introduced in
Definition \ref{def:piqi}.
\end{proposition}
\begin{proof} We begin by noting that any homomorphism in $\Hom (P\sb{n}, K\sb{m}\sp{1})$ can be regarded as a homomorphism in $\Hom (P\sb{n}, K\sb{m})$. On the other hand, the homomorphisms of $\Hom (P\sb{n}, K\sb{m})$ which cannot be regarded as homomorphisms from $\Hom (P\sb{n}, K\sb{m}\sp{1})$, are exactly those that map adjacent vertices of $P\sb{n}$ to the vertices of the exceptional edge. In this way,
\[
\hom (P\sb{n}, K\sb{m}\sp{1}) = \hom (P\sb{n}, K\sb{m}) - s\sb{n, m} = m(m-1)\sp{n-1} - s\sb{n, m} ,
\]
since $\hom (P\sb{n}, K\sb{m})$  is the chromatic polynomial of $P\sb{n}$,
\[
\hom (P\sb{n}, K\sb{m}) = \chi\sb{P\sb{n}} (m) = m(m-1)\sp{n-1}
\]
\cite[eq.\ (1.9), p.\ 3]{dkt}, and where $s\sb{n, m}$ is the cardinality of the set
\[
\{ f \in \Hom (P\sb{n}, K\sb{m}) \, | \, \text{ $f$
maps adjacent vertices to the vertices of the exceptional edge.} \}.
\]

Note further that $s\sb{n, m}=\sum\sb{k=1}\sp{n-1}s\sb{n, m}\sp{k}$, where $s\sb{n, m}\sp{k}$ is the cardinality of the set
\begin{multline*}
\{ f \in \Hom (P\sb{n}, K\sb{m}) \, | \\
\text{ $f$
maps exactly $k$ pairs of adjacent vertices to the vertices of the exceptional edge.} \}.
\end{multline*}
We will now prove that $s\sb{n, m}\sp{k}$ is given by the expression in the statement. We label vertices in $P\sb{n}$ $0, 1, \ldots , n-1$ in such a way that the edges are $01, 12, \ldots , (i-1)i, \ldots , (n-2)(n-1)$. Let $A$ and $B$ be the vertices of the exceptional edge in $K\sb{m}$. In order to specify that exactly $k$ pairs of adjacent vertices of $P\sb{n}$ are mapped to the pair $A, B$ in a given homomorphism, we specify the lowest vertex of each pair of these adjacent vertices. In this way, the sequence
\[
(0\leq) i\sb{1} < i\sb{2} < \cdots < i\sb{k} (\leq n-2)
\]
corresponds to the mapping of each pair of the $k$ adjacent vertices $i\sb{1}(i\sb{1}+1)$,  $i\sb{2}(i\sb{2}+1)$, ... , $i\sb{k}(i\sb{k}+1)$ injectively into the pair of adjacent vertices, $A$ and $B$. See Figure \ref{Fi:p8} for an illustrative example. It depicts the graph $P\sb{11}$ together with information turning it into a representative of  certain homomorphisms in $\Hom(P\sb{11}, K\sb{m})\setminus \Hom(P\sb{11}, K\sb{m}\sp{1})$. These homomorphisms map the vertices in full  to the exceptional vertices $A$, $B$. Specifically, vertices $0$, $1$; $3$, $4$, $5$; and $8$, $9$ are mapped in pairs injectively into $A, B$. The associated sequence of $i\sb{j}$'s is: $i\sb{1}=0, i\sb{2}=3, i\sb{3}=4, i\sb{4}=8$.

\begin{figure}[h!]
    \psfrag{1}{\huge $0=i\sb{1}$}
    \psfrag{2}{\huge $1$}
    \psfrag{3}{\huge $2$}
    \psfrag{4}{\huge $3=i\sb{2}$}
    \psfrag{5}{\huge $4=i\sb{3}$}
    \psfrag{6}{\huge $5$}
    \psfrag{7}{\huge $6$}
    \psfrag{8}{\huge $7$}
    \psfrag{9}{\huge $8=i\sb{4}$}
    \psfrag{10}{\huge $9$}
    \psfrag{11}{\huge $10$}
    \psfrag{...}{\huge $\ldots $}
    \centerline{\scalebox{.50}{\includegraphics{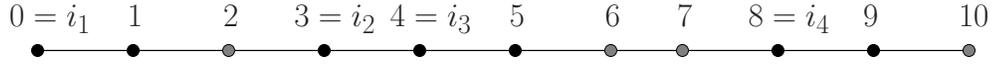}}}
    \caption{Graph $P\sb{11}$ together with an assignment to the vertices of the exceptional edge.}\label{Fi:p8}
\end{figure}

Note that $i\sb{j}-i\sb{j-1}\neq 2$ for $j=2, \ldots , k$. In fact, if $i\sb{j\sb{0}}-i\sb{j\sb{0}-1} = 2$, then assume $i\sb{j\sb{0}-1}$ and $i\sb{j\sb{0}}$ are mapped to the same element of $\{ A, B  \}$. Then, $i\sb{j\sb{0}-1}+1$ is mapped to the other element of $\{ A, B  \}$, in this way contributing with another pair of adjacent vertices, $i\sb{j\sb{0}-1}+1, i\sb{j\sb{0}}$, being mapped to $A, B$,
thereby increasing the number of the pairs of these vertices to $k+1$, which contradicts the assumption. Assume now $i\sb{j\sb{0}-1}$ and $i\sb{j\sb{0}}$ are mapped to distinct elements of $\{ A, B  \}$. In this case, $i\sb{j\sb{0}-1}+1$ has to be mapped to one of them, thereby contradicting preservation of adjacency.

Given a sequence of $k$ pairs being mapped to vertices $A$ and $B$ via their sequence of lower vertices, say
\[
(0\leq) i\sb{1} < i\sb{2} < \cdots < i\sb{k} (\leq n-2),
\]
and $i\sb{j}-i\sb{j-1}\neq 2$ for $j=2, \ldots , k$, we call \emph{cluster} each maximal subsequence of the latter sequence with respect to the property that each two consecutive terms of the subsequence differ by $1$ , i.e., $i\sb{j\sb{r}}-i\sb{j\sb{r}-1}=1$, together with the vertex next to the last term of the maximal subsequence at issue.  That is, if $i\sb{j\sb{0}}, i\sb{j\sb{0}}+1, i\sb{j\sb{0}}+2, \ldots , i\sb{j\sb{0}}+m-1$ is such a maximal sequence, then the cluster it gives rise to is $i\sb{j\sb{0}}, i\sb{j\sb{0}}+1, i\sb{j\sb{0}}+2, \ldots , i\sb{j\sb{0}}+m-1 , i\sb{j\sb{0}}+m$. In the situation depicted in Figure \ref{Fi:p8} there are $3$ maximal sequences: $0$; $3, 4$; and $8$. These give rise to the clusters $0, 1$; $3, 4, 5$; and $8, 9$, respectively.

There are two possible ways of mapping the elements of a cluster, say $i\sb{j\sb{0}}, i\sb{j\sb{0}}+1, i\sb{j\sb{0}}+2, \ldots , i\sb{j\sb{0}}+m$, to $A$ and $B$. Either
\[
i\sb{j\sb{0}} \mapsto A,\qquad (i\sb{j\sb{0}}+1) \mapsto B, \qquad (i\sb{j\sb{0}}+2) \mapsto A,\qquad (i\sb{j\sb{0}}+3) \mapsto B, \quad \ldots
\]

or

\[
i\sb{j\sb{0}} \mapsto B,\qquad (i\sb{j\sb{0}}+1) \mapsto A, \qquad (i\sb{j\sb{0}}+2) \mapsto B,\qquad (i\sb{j\sb{0}}+3) \mapsto A, \quad \ldots
\]

Now, assume $i\sb{j\sb{0}}, i\sb{j\sb{0}}+1, i\sb{j\sb{0}}+2, \ldots , i\sb{j\sb{0}}+m$ and $i\sb{j'\sb{0}}, i\sb{j'\sb{0}}+1, i\sb{j'\sb{0}}+2, \ldots , i\sb{j'\sb{0}}+m'$ are two consecutive clusters with $i\sb{j\sb{0}}+m < i\sb{j'\sb{0}}$. We call \emph{gap} the sequence $i\sb{j\sb{0}}+m+1 < i\sb{j\sb{0}}+m+2 < \cdots < i\sb{j'\sb{0}}-1$ that lies between consecutive clusters. Should $i\sb{1}>0$, we also call \emph{gap} the sequence of vertices right before the first cluster, $0, \ldots , i\sb{1}-1$ . Should $i\sb{k}<n-2$, we also call \emph{gap} the sequence of vertices right after the last cluster, $i\sb{k}+2, \ldots , n-1$. In the situation depicted in Figure \ref{Fi:p8} there are three gaps: $2$; $6, 7$; and $10$.

Then, for each homomorphism from $P\sb{n}$ to $K\sb{m}$ that maps exactly $k$ pairs of adjacent vertices into the vertices $A$ and $B$, there is a sequence of lower vertices as explained above. This sequence gives rise to clusters and gaps. Specifically, the sequence of lower vertices induces a partition of $P\sb{n}$ into clusters and gaps. We already saw that there are two ways of mapping the vertices in a cluster to $A$ and $B$. Then, for each cluster, there is a factor of two in the enumeration of these homomorphisms. We now analyze the contribution of the gaps to the counting.

Suppose a gap only has one vertex. Then, this vertex cannot be mapped either to $A$ or to $B$, for otherwise it would either increase the number of adjacent vertices being mapped to $A$ and $B$, or it would contradict the preservation of adjacency. This vertex can, then, be mapped to one in $m-2$ vertices in $K\sb{m}$.

We, now, consider gaps nested between two consecutive clusters. Later, we will consider gaps which make contact with only one cluster: $0, 1, \ldots , i\sb{1}-1$, when $i\sb{1}>0$, and/or $i\sb{k}+2, \ldots , n-1$, when $i\sb{k}<n-2$.

Suppose a gap has two elements, say $g\sb{1} < g\sb{2}$. Then, $g\sb{1}$ cannot be mapped to either $A$ or $B$, for otherwise it would be included in the cluster preceding it. So, $g_1$ can be mapped to one in $m-2$ vertices. Then, $g\sb{2}$ cannot be mapped to the vertex $g\sb{1}$ has been mapped to, and it also cannot be mapped to $A$ or to $B$. But it can  be mapped to one in $m-3$ vertices. The contribution of a two element gap is, then, $(m-2)(m-3)$.

Suppose a gap has three elements, say $g\sb{1} < g\sb{2} < g\sb{3}$. As discussed before, $g\sb{1}$ and $g\sb{3}$ cannot be mapped to $A$ or to $B$. But $g\sb{2}$ can be mapped to $A$ or to $B$. Note, however, that if $g\sb{2}$ is not mapped either to $A$ or to $B$ then $g\sb{3}$ cannot be mapped also to the vertex $g\sb{2}$ has been mapped to, whereas if $g\sb{2}$ has been mapped to $A$ or to $B$ then $g\sb{3}$ can be mapped to all vertices but $A$ or $B$. In this way the contribution for the number of homomorphisms from a gap with three elements is $(m-2)2(m-2)+(m-2)(m-3)(m-3)$.

Suppose a gap has four elements, $g\sb{1} < g\sb{2} < g\sb{3} < g\sb{4}$. We believe this will illustrate what happens in the general case. As above, $g\sb{1}$ and $g\sb{4}$ cannot be mapped to $A$ or to $B$. If $g\sb{2}$ is mapped to $A$ or to $B$ ($2$ possibilities), then $g\sb{3}$ can be mapped to any but $A$ or $B$ ($m-2$ possibilities) and, then, $g\sb{4}$ cannot be mapped to the vertex $g\sb{3}$ has been mapped to, nor to $A$ nor to $B$ ($m-3$ possibilities). In this case, we have $(m-2)2(m-2)(m-3)$. If $g\sb{2}$ is not mapped to $A$ nor to $B$, nor to the vertex $g\sb{1}$ has been mapped to ($m-3$ possibilities), then $g\sb{3}$ can be mapped to $A$ or to $B$ ($2$ possibilities), in which case $g\sb{4}$ can be mapped to any but $A$ or $B$ ($m-2$ possibilities). In this case, we have $(m-2)(m-3)2(m-2)$. Finally, if $g\sb{2}$ is not mapped to $A$ nor to $B$, nor to the vertex $g\sb{1}$ has been mapped to ($m-3$ possibilities), and $g\sb{3}$ is not mapped to $A$ nor $B$, nor the vertex $g\sb{2}$ has been mapped to ($m-3$ possibilities), then $g\sb{4}$ has only $m-3$ possibilities. In this case, we have $(m-2)(m-3)(m-3)(m-3)$ possibilities. The contribution for the number of homomorphisms from a gap with four elements is, then, $(m-2)2(m-2)(m-3)+(m-2)(m-3)2(m-2)+(m-2)(m-3)(m-3)(m-3)$.

\begin{figure}[h!]
    \psfrag{0}{\huge $\mathbf 0$}
    \psfrag{1}{\huge $\begin{cases}\mathbf 0\\ \mathbf 1\end{cases}$}
    \psfrag{ij-1}{\huge $i\sb{j-1}$}
    \psfrag{ij-1+1}{\huge $i\sb{j-1}+1$}
    \psfrag{ij}{\huge $i\sb{j}$}
    \psfrag{r1}{\huge $g\sb{1}$}
    \psfrag{r2}{\huge $g\sb{2}$}
    \psfrag{r3}{\huge $g\sb{3}$}
    \psfrag{r4}{\huge $g\sb{4}$}
    \psfrag{rl-1}{\huge $g\sb{l-1}$}
    \psfrag{rl}{\huge $g\sb{l}$}
    \psfrag{...}{\huge $\ldots $}
    \centerline{\scalebox{.50}{\includegraphics{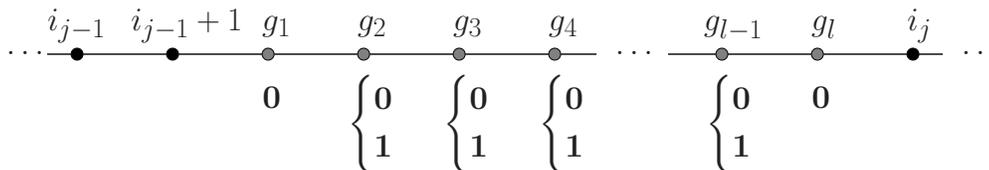}}}
    \caption{$l$ vertices between $i\sb{j-1}+1$ and $i\sb{j}$.}\label{Fi:prooftree}
\end{figure}

Figure \ref{Fi:prooftree} now illustrates the general case. We distinguish two sorts of assignments for each vertex of a gap (designated $g\sb{j}'s$). Either the vertex is assigned to one of $A$ or $B$ (the $``1"$ assignment), or the vertex is not assigned to
$A$ or to $B$ (the $``0"$ assignment). Then $0$'s are always associated to the first and last vertices of a gap, $g\sb{1}$ and $g\sb{l}$. The other vertices of the gap can, however, have both $0$ or $1$ assignments. Let us now see how these assignments contribute to the counting of homomorphisms that map exactly $k$ pairs of adjacent vertices to $A$ and $B$. The first vertex of the gap can be mapped to one of $m-2$ vertices since it cannot be mapped to $A$ or to $B$. Then, $g\sb{2}$ either receives a type $1$ assignment ($2$ possibilities), or a type $0$ assignment ($m-3$ possibilities). Then, the possibilities for $g\sb{3}$ are as follows. If $g\sb{2}$ received a $1$ assignment, then $g\sb{3}$ has to receive a $0$ assignment with $m-2$ possibilities. If $g\sb{2}$ received a $0$ assignment, then $g\sb{3}$ can receive either one of the two types of assignments: if $g\sb{3}$ receives a type $1$ assignment there are $2$ possibilities; if $g\sb{3}$ receives a type $0$ assignment there are $m-3$ possibilities. So from $1$'s stem $0$'s with $m-2$ possibilities. From $0$'s stem either $0$'s with $m-3$ possibilities or $1$'s with $2$ possibilities. Note that the product of these factors gives the summands of the polynomial $p\sb{l}(m)$ (see Definition \ref{def:piqi}) from the graph $G\sb{l}$ (see Definition \ref{def:graphs}). More precisely, each choice of $0$ or $1$ at each vertex of the gap, subject to the first and the last being assigned $0$'s, yields a maximal path in the graph $G\sb{l}$ which gives rise to the corresponding summand of the $p\sb{l}(m)$ polynomial. Then $p\sb{l}(m)$ yields all the possible contributions of the gap with $l$ vertices to the homomorphisms being enumerated.

Finally, recall that each homomorphism with exactly $k$ pairs of adjacent vertices being mapped to the adjacent vertices $A$ and $B$ corresponds to a sequence of lower vertices
\[
(0\leq) i\sb{1} < i\sb{2} < \cdots < i\sb{k} (\leq n-2),
\]
with $i\sb{j}-i\sb{j-1}\neq 2$ for $j=2, \ldots , k$. Then each gap nested between two clusters corresponds to a $j$ such that  $i\sb{j}-i\sb{j-1} > 2$ for some $j=2, \ldots , k$. The $p\sb{l}(m)$ polynomial corresponding to this gap has $l = i\sb{j} - i\sb{j-1} - 2$.

Notwithstanding, the gaps may not be nested between two consecutive clusters. Consider the case $i\sb{k}<n-2$. There will be a gap of $n-2-i\sb{k}$ vertices with a cluster to its left but no cluster to its right. For this particular gap, the counting of possibilities is yielded by $q\sb{n-2-i\sb{k}}(m)$, where the $q\sb{l}(m)$'s are the polynomials introduced in Definition \ref{def:piqi}. In fact, the attribution of possibilities to the vertices of this gap is the same as for the other gaps, except for the last vertex of this gap, vertex $n-1$. Since there is no cluster next to this vertex, then the possibilities attributed to it are either $m-1$, if the preceding vertex has had a $0$ assignment, or $m-2$, if the preceding vertex has had a $1$ assignment. An analogous situation occurs for $i\sb{1}>0$. In this case, the counting of the possibilities is yielded by $q_{i_1}(m)$. In the expression for $s\sb{n, m}\sp{k}$ below, we single out the instances where these latter gaps occur.

Then,
\[
s\sb{n, m}\sp{n-1}  = 2
\]
since in this $k=n-1$ situation there is only one cluster (which involves all the vertices of the path $P\sb{n}$) and there is no gap. For $1\leq k \leq n-2$,
\begin{align*}
s\sb{n, m}\sp{k} & = \frac{1}{2}\sum\sb{\substack{\{ 0 = i\sb{1}<i\sb{2}< \cdots < i\sb{k} = n-2 \} \\
\text{ s.t. } i\sb{j}-i\sb{j-1}\neq 2}}
 \Biggl[  \prod\sb{\substack{i\sb{j} \, \text{ s.t. }\\
i\sb{j}-i\sb{j-1}>2 \\  2\leq j \leq k}} 2\,
p\sb{i\sb{j}-i\sb{j-1}-2}(m) \Biggr]   \\
& \qquad + \sum\sb{\substack{\{ 0 = i\sb{1}<i\sb{2}< \cdots < i\sb{k} < n-2 \} \\
\text{ s.t. } i\sb{j}-i\sb{j-1}\neq 2}}
\Biggl[  \prod\sb{\substack{i\sb{j} \, \text{ s.t. }\\
i\sb{j}-i\sb{j-1}>2 \\  2\leq j \leq k}} 2\,
p\sb{i\sb{j}-i\sb{j-1}-2}(m) \Biggr]  \cdot q\sb{n-2-i\sb{k}}(m)  \\
& \qquad + \sum\sb{\substack{\{ 0 < i\sb{1}<i\sb{2}< \cdots < i\sb{k} = n-2 \} \\
\text{ s.t. } i\sb{j}-i\sb{j-1}\neq 2}}q\sb{i\sb{1}}(m)
\cdot \Biggl[  \prod\sb{\substack{i\sb{j} \, \text{ s.t. }\\
i\sb{j}-i\sb{j-1}>2 \\  2\leq j \leq k}} 2\,
p\sb{i\sb{j}-i\sb{j-1}-2}(m) \Biggr]    \\
& \qquad + 2 \sum\sb{\substack{\{ 0 < i\sb{1}<i\sb{2}< \cdots < i\sb{k} < n-2 \} \\
\text{ s.t. } i\sb{j}-i\sb{j-1}\neq 2}}q\sb{i\sb{1}}(m)
\cdot \Biggl[  \prod\sb{\substack{i\sb{j} \, \text{ s.t. }\\
i\sb{j}-i\sb{j-1}>2 \\  2\leq j \leq k}} 2\,
p\sb{i\sb{j}-i\sb{j-1}-2}(m) \Biggr]  \cdot q\sb{n-2-i\sb{k}}(m),
\end{align*}
and the $p\sb{i}$'s and $q\sb{i}$'s are the polynomials introduced in
Definition \ref{def:piqi}. A $2$ right before a $p$ polynomial stands for the two possibilities for the cluster right before (or right after) the gap the $p$ polynomial stands for. The $\frac{1}{2}$ in the first summand above and the $2$ factor before the last summand above stand for corrections to these enumerations of gaps when, respectively, each gap is nested between two consecutive clusters, and when there are exactly two gaps which are not nested between two consecutive clusters.

This completes the proof.

\end{proof}

\subsection{The source graphs are cycles}\label{subsect:cycle}

We remark that a \emph{cycle on $n$ vertices}, denoted $C\sb{n}$, is a graph with vertices $0, 1, \ldots , n-1$ and edges $01, 12, \ldots , (n-2)(n-1), (n-1)0$. We take $n \geq 4$ for $C\sb{3}$ is isomorphic with $K\sb{3}$.

\begin{proposition}\label{prop:cnkm1}For integers $n\geq 4$, $m\geq 3$,
\[
\hom(C\sb{n}, K\sb{m}\sp{1}) = (m-1)\biggl( (m-1)\sp{n-1}+(-1)\sp{n}\biggr)  -
\sum\sb{k=1}\sp{n}t\sb{n, m}\sp{k},
\]
where
\begin{align*}
t\sb{n, m}\sp{n}  &= 2,  \qquad \text{ if $n$ is even},\\
t\sb{n, m}\sp{n}  &= 0,  \qquad \text{ if $n$ is odd},\\
t\sb{n, m}\sp{n-1} &= 0,
\end{align*}
and, for $1\leq k \leq n-2$,
\begin{align*}
t\sb{n, m}\sp{k} & = \sum\sb{\substack{\{ 0 = i\sb{1}<i\sb{2}< \cdots < i\sb{k} = n-1 \} \\
\text{ s.t. } i\sb{j}-i\sb{j-1}\neq 2}}
 \Biggl[  \prod\sb{\substack{i\sb{j} \, \text{ s.t. }\\
i\sb{j}-i\sb{j-1}>2 \\  2\leq j \leq k}} 2\,
p\sb{i\sb{j}-i\sb{j-1}-2}(m) \Biggr]  \\
& \qquad + \sum\sb{\substack{\{ 0 = i\sb{1}<i\sb{2}< \cdots < i\sb{k} < n-2 \} \\
\text{ s.t. } i\sb{j}-i\sb{j-1}\neq 2}}
\Biggl[  \prod\sb{\substack{i\sb{j} \, \text{ s.t. }\\
i\sb{j}-i\sb{j-1}>2 \\  2\leq j \leq k}} 2\,
p\sb{i\sb{j}-i\sb{j-1}-2}(m) \Biggr]  \cdot 2p\sb{n-i\sb{k}-2}(m)  \\
& \qquad + \sum\sb{\substack{\{ 1 < i\sb{1}<i\sb{2}< \cdots < i\sb{k} = n-1 \} \\
\text{ s.t. } i\sb{j}-i\sb{j-1}\neq 2}}2p\sb{i\sb{1}}(m)
\cdot \Biggl[  \prod\sb{\substack{i\sb{j} \, \text{ s.t. }\\
i\sb{j}-i\sb{j-1}>2 \\  2\leq j \leq k}} 2\,
p\sb{i\sb{j}-i\sb{j-1}-2}(m) \Biggr]    \\
& \qquad +  \sum\sb{\substack{\{ 0 < i\sb{1}<i\sb{2}< \cdots < i\sb{k} < n-1 \} \\
\text{ s.t. } i\sb{j}-i\sb{j-1}\neq 2}}2p\sb{i\sb{1}+n-i\sb{k}-2}(m)
\cdot \Biggl[  \prod\sb{\substack{i\sb{j} \, \text{ s.t. }\\
i\sb{j}-i\sb{j-1}>2 \\  2\leq j \leq k}} 2\,
p\sb{i\sb{j}-i\sb{j-1}-2}(m) \Biggr] ,
\end{align*}
and the $p\sb{i}$'s are the polynomials introduced in
Definition \ref{def:piqi}.
\end{proposition}
\begin{proof} The proof is similar to that of Proposition \ref{prop:pnkm1} and will be largely omitted.

We remark that
\[
\Hom(C\sb{n}, K\sb{m}) = \chi_{C\sb{n}}(m) = (m-1)\biggl( (m-1)\sp{n-1}+(-1)\sp{n}\biggr)
\]
\cite[eq.\ (1.13), p.\ 7]{dkt}.

Observing that $\Hom(C\sb{n}, K\sb{m}\sp{1})$ embeds in $\Hom(C\sb{n}, K\sb{m})$, we just have to enumerate the homomorphisms from the latter set which are not from the former set; we denote these $t\sb{n, m}$. As before, we obtain $t\sb{n, m}$ by counting those homomorphisms that map exactly $k$ pairs of adjacent vertices in $C\sb{n}$ to $K\sb{m}$, and, then, adding over $k$, $t\sb{n, m}=\sum\sb{k=1}\sp{n}t\sb{n, m}\sp{k}$.

The notion of cluster, here, has to be somewhat specialized to encompass the cases of $k$ pairs of vertices with $i\sb{1}=0$ and $i\sb{k}=n-1$. In these situations, the maximal subsequence starting at $i\sb{1}=0$ merges with the maximal subsequence ending at $i\sb{k}=n-1$, forming only one cluster. Especially, when there is only one cluster and there are vertices of the cycle, $C\sb{n}$, which are not part of the cluster, then we call gap the set of these vertices. Note that this gap is nested inside the only cluster.

If $n$ is even, then there is a cluster without gaps for $k=n$; this cluster involves all the vertices in the cycle $C\sb{n}$. If $n$ is odd, there is no such situation for otherwise adjacency would be violated. If $k = n-1$, then $t\sb{n, m}\sp{k}=0$ for the pair of vertices apparently left out would contribute another pair of vertices or violate adjacency.

 For $1\leq k < n-1$, each gap is nested between two clusters or nested inside one cluster (when there is only one cluster and only one gap).

 Concerning the homomorphisms in $\Hom (C\sb{n}, K\sb{m}) \setminus \Hom (C\sb{n}, K\sb{m}\sp{1})$, clusters can be mapped in one of two ways and the mappings of gaps are ruled by the graphs introduced in Definition \ref{def:graphs} in conjunction with the polynomials $p\sb{l}(m)$ introduced in Definition \ref{def:piqi} ($l$ is the number of the vertices in the gap at issue).

 We leave the remaining details to the reader.

\end{proof}

\begin{corollary}\label{cor:coddk31} For integer $n'\geq 2$,
\[
\hom(C\sb{2n'+1}, K\sb{3}\sp{1}) = 0.
\]
\end{corollary}
\begin{proof} The homomorphic image of a cycle on an odd number of vertices is a cycle \cite{hell}. Since there is no cycle in $K\sb{3}\sp{1}$, this concludes the proof.
\end{proof}

\subsection{The source graphs are broken wheels}\label{subsect:bwheel}

Let $n$ be an integer, $n\geq 3$. We remark that a \emph{broken wheel on $n$ spokes}, denoted $BW\sb{n}$, is a graph with $n+1$ vertices, formed by a path on $n$ vertices plus an extra vertex, called the hub, which is adjacent to every other vertex. The hub is denoted $0$ and the other vertices $1, 2, \ldots , n$; the edges are $01, 02, 03, \ldots , 0n$ (the spokes) and $12, 23, 34, \ldots , (n-1)n$.

\begin{proposition}\label{prop:bwnkm1}For integers $n\geq 3$, $m\geq 3$,
\[
\hom(BW\sb{n}, K\sb{m}\sp{1}) = m(m-1)(m-2)\sp{n-1} -
\sum\sb{k=1}\sp{n-1}u\sb{n, m}\sp{k},
\]
where
\[
u\sb{n, m}\sp{k}  = u\sb{n, m}\sp{k, 0}+u\sb{n, m}\sp{k, 1},
\]
and, keeping the notation of Proposition \ref{prop:pnkm1},
\[
u\sb{n, m}\sp{k, 0}  = (m-2)s\sb{n, m-1}\sp{k} \qquad  1 \leq k \leq n-1,
\]
and, for $1\leq k \leq \lceil n/2 \rceil$,
\begin{align*}
&u\sb{n, m}\sp{k, 1}   = 2\sum\sb{\substack{\{  1= i\sb{1} < i\sb{2} < i\sb{3} < \cdots < i\sb{k} = n  \} \\
\text{ s.t. } i\sb{j}-i\sb{j-1}>1, \;\;  2\leq j\leq k }}
 \Biggl[  \prod\sb{\substack{i\sb{j} \, \text{ s.t. }\\
i\sb{j}-i\sb{j-1}(>1) \\  2\leq j \leq k}} (m-2)(m-3)\sp{i\sb{j}-i\sb{j-1}-2} \Biggr]   \\
& \quad + 2\sum\sb{\substack{\{  1 < i\sb{1} < i\sb{2} < i\sb{3} < \cdots < i\sb{k} = n  \} \\
\text{ s.t. } i\sb{j}-i\sb{j-1}>1, \;\;  2\leq j\leq k }}
 (m-2)(m-3)\sp{i\sb{1}-2}\cdot \Biggl[  \prod\sb{\substack{i\sb{j} \, \text{ s.t. }\\
i\sb{j}-i\sb{j-1}(>1) \\  2\leq j \leq k}} (m-2)(m-3)\sp{i\sb{j}-i\sb{j-1}-2} \Biggr]  \\
& \quad + 2\sum\sb{\substack{\{  1 = i\sb{1} < i\sb{2} < i\sb{3} < \cdots < i\sb{k} < n  \} \\
\text{ s.t. } i\sb{j}-i\sb{j-1}>1, \;\;  2\leq j\leq k }}
 \Biggl[  \prod\sb{\substack{i\sb{j} \, \text{ s.t. }\\
i\sb{j}-i\sb{j-1}(>1) \\  2\leq j \leq k}} (m-2)(m-3)\sp{i\sb{j}-i\sb{j-1}-2} \Biggr] \cdot  (m-2)(m-3)\sp{n-i\sb{k}-1}  \\
& \quad + 2\sum\sb{\substack{\{  1 < i\sb{1} < i\sb{2} < i\sb{3} < \cdots < i\sb{k} < n  \} \\
\text{ s.t. } i\sb{j}-i\sb{j-1}>1, \;\;  2\leq j\leq k }}
 (m-2)\sp{2}(m-3)\sp{i\sb{1}+n-i\sb{k}-3}\cdot \Biggl[  \prod\sb{\substack{i\sb{j} \, \text{ s.t. }\\
i\sb{j}-i\sb{j-1}(>1) \\  2\leq j \leq k}} (m-2)(m-3)\sp{i\sb{j}-i\sb{j-1}-2} \Biggr] ,
\end{align*}
whereas
\[
u\sb{n, m}\sp{k, 1}=0 \qquad \text{ for } \lceil n/2 \rceil < k \leq n .
\]
\end{proposition}
\begin{proof}

We remark that
\[
\Hom(BW\sb{n}, K\sb{m}) = \chi_{BW\sb{n}}(m) = m(m-1)(m-2)\sp{n-1} ,
\]
by applying the fundamental recursion theorem \cite[eq.\ (1.11), p.\ 4]{dkt} to the chromatic polynomial for the wheel graphs \cite[eq.\ (1.23), p.\ 13]{dkt}.

The strategy will be to observe that $\Hom (BW\sb{n}, K\sb{m}\sp{1})$ embeds in $\Hom (BW\sb{n}, K\sb{m})$ and then to count which homomorphisms of the latter set are not in the former set; these are the ones that map at least a pair of adjacent vertices of $BW\sb{n}$ to the exceptional edge $AB$ in $K\sb{m}$. We count these by counting how many homomorphisms there are which map exactly $k$ pairs of adjacent vertices of $BW\sb{n}$ to $K\sb{m}$ (we denote these numbers $u\sb{n, m}\sp{k}$), and, then, adding over $k$ from $1$ to $n-1$.

Concerning these homomorphisms which map adjacent vertices of $BW\sb{n}$ into the vertices of the exceptional edge of $K\sb{m}$, we distinguish two cases. In one case, the hub is mapped to one of the vertices of the exceptional edge, and, in the other one, the hub is not mapped to any of these vertices.

\begin{figure}[h!]
    \psfrag{l}{\huge $\mathbf 0$}
    \psfrag{j}{\huge $\begin{cases}\mathbf 0\\ \mathbf 1\end{cases}$}
    \psfrag{0}{\huge $0$}
    \psfrag{1}{\huge $1$}
    \psfrag{2}{\huge $2$}
    \psfrag{3}{\huge $3$}
    \psfrag{4}{\huge $4$}
    \psfrag{5}{\huge $5$}
    \psfrag{6}{\huge $6$}
    \psfrag{rl-1}{\huge $g\sb{l-1}$}
    \psfrag{rl}{\huge $g\sb{l}$}
    \psfrag{...}{\huge $\ldots $}
    \centerline{\scalebox{.50}{\includegraphics{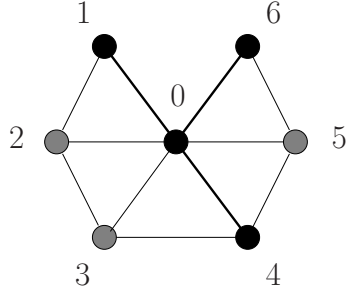}}}
    \caption{The broken wheel on $6$ spokes, $BW\sb{6}$. The vertices drawn in full are assigned to the vertices of the exceptional edge.}\label{Fi:bw6}
\end{figure}

 Suppose that, for a given $f\in \Hom (BW\sb{n}, K\sb{m})$, the hub is part of a pair of vertices in $BW\sb{n}$ which is mapped to the exceptional edge. Then each pair of adjacent vertices of $BW\sb{n}$ in this homomorphism which is mapped to the $A$, $B$ vertices, involves the hub. Otherwise, adjacency would not be preserved, since every other vertex is adjacent to the hub. In short, if one spoke is used, then only spokes are used. This is the $1$ situation, as signaled by the upper index $1$ in $u\sb{n, m}\sp{k, 1}$. Moreover, no two consecutive vertices in the path part of the $BW\sb{n}$ can be mapped to either of the vertices of the exceptional edge for otherwise, adjacency would not be preserved. Consider one such homomorphism which maps $k$ pairs of vertices to the vertices of the exceptional edge. This homomorphism is clearly identified by describing the vertices on the path part of the broken wheel that are mapped to the other vertex of the pair of exceptional vertices and conversely. On the other hand, each such set of vertices is in one-to-one correspondence with the subsets of $k$ elements from $\{ 1, 2, \ldots , n  \}$ (the path part of the broken wheel) such that any two elements are more than one unit apart. This sets boundaries on the values of $k$: $1\leq k \leq \lceil n/2 \rceil$. The rest of the reasoning is analogous (although not the same) to the reasoning associated with (clusters and) gaps made above, in the proofs of Propositions \ref{prop:pnkm1} and \ref{prop:cnkm1}. In Figure \ref{Fi:bw6}, we display the broken wheel on $6$ spokes along with a representative of a homomorphism to $K\sb{m}$ mapping $3$ pairs of adjacent vertices to the vertices of the exceptional edge and each of these pairs involve the hub, vertex $0$. The vertices that are mapped to the vertices of the exceptional edge are drawn in full, and so are the edges connecting them. With the notation used in the expressions for $u\sb{n, m}\sp{k, 1}$, this corresponds to $i\sb{1}=1, i\sb{2}=4, i\sb{3}=6$. Thus, either the hub is mapped to $A$ and vertices $1$, $4$, and $6$ are mapped to $B$, or the hub is mapped to $B$, and vertices $1$, $4$, and $6$ are mapped to $A$, accounting so far for two possibilities. Then, vertex $2$ cannot be mapped to $A$ nor to $B$ ($m-2$ possibilities); vertex $3$ cannot be mapped to $A$, nor to $B$, nor to the vertex $2$ has been mapped to ($m-3$ possibilities). Finally, vertex $5$ cannot be mapped to $A$ nor to $B$ ($m-2$ possibilities). In this way, Figure \ref{Fi:bw6} stands for $2\cdot (m-2)(m-3)\cdot (m-2)$ homomorphisms which map exactly $3$ pairs of adjacent vertices, each one of which involving the hub, to the vertices of the exceptional edge.

 Then, for $1\leq k \leq \lceil n/2 \rceil$,
\begin{align*}
&u\sb{n, m}\sp{k, 1}   = 2\sum\sb{\substack{\{  1= i\sb{1} < i\sb{2} < i\sb{3} < \cdots < i\sb{k} = n  \} \\
\text{ s.t. } i\sb{j}-i\sb{j-1}>1, \;\;  2\leq j\leq k }}
 \Biggl[  \prod\sb{\substack{i\sb{j} \, \text{ s.t. }\\
i\sb{j}-i\sb{j-1}(>1) \\  2\leq j \leq k}} (m-2)(m-3)\sp{i\sb{j}-i\sb{j-1}-2} \Biggr]   \\
& \quad + 2\sum\sb{\substack{\{  1 < i\sb{1} < i\sb{2} < i\sb{3} < \cdots < i\sb{k} = n  \} \\
\text{ s.t. } i\sb{j}-i\sb{j-1}>1, \;\;  2\leq j\leq k }}
 (m-2)(m-3)\sp{i\sb{1}-2}\cdot \Biggl[  \prod\sb{\substack{i\sb{j} \, \text{ s.t. }\\
i\sb{j}-i\sb{j-1}(>1) \\  2\leq j \leq k}} (m-2)(m-3)\sp{i\sb{j}-i\sb{j-1}-2} \Biggr]   \\
& \quad + 2\sum\sb{\substack{\{  1 = i\sb{1} < i\sb{2} < i\sb{3} < \cdots < i\sb{k} < n  \} \\
\text{ s.t. } i\sb{j}-i\sb{j-1}>1, \;\;  2\leq j\leq k }}
 \Biggl[  \prod\sb{\substack{i\sb{j} \, \text{ s.t. }\\
i\sb{j}-i\sb{j-1}(>1) \\  2\leq j \leq k}} (m-2)(m-3)\sp{i\sb{j}-i\sb{j-1}-2} \Biggr] \cdot  (m-2)(m-3)\sp{n-i\sb{k}-1}  \\
& \quad + 2\sum\sb{\substack{\{  1 < i\sb{1} < i\sb{2} < i\sb{3} < \cdots < i\sb{k} < n  \} \\
\text{ s.t. } i\sb{j}-i\sb{j-1}>1, \;\;  2\leq j\leq k }}
 (m-2)\sp{2}(m-3)\sp{i\sb{1}+n-i\sb{k}-3}\cdot \Biggl[  \prod\sb{\substack{i\sb{j} \, \text{ s.t. }\\
i\sb{j}-i\sb{j-1}(>1) \\  2\leq j \leq k}} (m-2)(m-3)\sp{i\sb{j}-i\sb{j-1}-2} \Biggr] ,
\end{align*}
where the $2$ factors stand for the two possibilities of mapping the hub to, $A$ or $B$.

Suppose, now, that, for a given $f\in \Hom (BW\sb{n}, K\sb{m})$, the hub is not part of any pair of vertices which are mapped to the exceptional edge. This is the $0$ situation, as signaled by the upper index $0$ in $u\sb{n, m}\sp{k, 0}$. Then, only vertices on the path part of the broken wheel are mapped to the vertices of the exceptional edge. Since the hub can be mapped to every vertex in $K\sb{m}$ but $A$ or $B$, there are $m-2$ possibilities for it. Let us, now, calculate the possibilities for exactly $k$ pairs of adjacent vertices on the path part of the broken wheel to be mapped to the vertices in the exceptional edge. Since the hub has already taken up one of the vertices we are left with $m-1$ vertices in the target, otherwise violating preservation of adjacency. The possibilities for the path part of the broken wheel are, then, $s\sb{n, m-1}\sp{k}$, using the notation of Proposition \ref{prop:pnkm1}. The total number of possibilities in this $0$ situation of mapping $k$ adjacent vertices to $A$ and to $B$ is, then,
\[
u\sb{n, m}\sp{k, 0} = (m-2)\cdot s\sb{n, m-1}\sp{k} .
\]
This completes the proof.
\end{proof}

\begin{corollary}\label{cor:bwoddk31} For integer $n'\geq 1$,
\[
\hom(BW\sb{2n'+1}, K\sb{3}\sp{1}) = 0 .
\]
\end{corollary}
\begin{proof} Omitted.

\end{proof}

\subsection{The source graphs are wheels}\label{subsect:wheel}

Let $n$ be an integer, $n\geq 3$. We remark that a \emph{wheel on $n$ spokes}, denoted $W\sb{n}$, is a graph consisting of $n+1$ vertices formed by a cycle on $n$ vertices plus an extra vertex, called the hub, which is adjacent to every other vertex. The hub is denoted $0$, and the other vertices $1, 2, \ldots , n$; the vertices are $01, 02, 03, \ldots , 0n$ (the spokes) and $12, 23, 34, \ldots , (n-1)n, n1$.

\begin{proposition}\label{prop:wnkm1}For integers $n\geq 3$, $m\geq 3$,
\[
\hom(W\sb{n}, K\sb{m}\sp{1}) = m(m-2)\big[ (m-2)\sp{n-1}+(-1)\sp{n}\big] -
\sum\sb{k=1}\sp{n}v\sb{n, m}\sp{k} ,
\]
where
\[
v\sb{n, m}\sp{k}  = v\sb{n, m}\sp{k, 0}+v\sb{n, m}\sp{k, 1} ,
\]
and, keeping the notation of Proposition \ref{prop:cnkm1},
\[
v\sb{n, m}\sp{k, 0}  = (m-2)t\sb{n, m-1}\sp{k}, \qquad 1 \leq k \leq n ,
\]
and, for $1 \leq k \leq \lfloor n/2 \rfloor$,
\begin{align*}
&v\sb{n, m}\sp{k, 1}   =  2\sum\sb{\substack{\{  1 < i\sb{1} < i\sb{2} < i\sb{3} < \cdots < i\sb{k} = n  \} \\
\text{ s.t. } i\sb{j}-i\sb{j-1}>1, \;\;  2\leq j\leq k }}
 (m-2)(m-3)\sp{i\sb{1}-2}\cdot \Biggl[  \prod\sb{\substack{i\sb{j} \, \text{ s.t. }\\
i\sb{j}-i\sb{j-1}(>1) \\  2\leq j \leq k}} (m-2)(m-3)\sp{i\sb{j}-i\sb{j-1}-2} \Biggr]   \\
& \quad + 2\sum\sb{\substack{\{  1 = i\sb{1} < i\sb{2} < i\sb{3} < \cdots < i\sb{k} < n  \} \\
\text{ s.t. } i\sb{j}-i\sb{j-1}>1, \;\;  2\leq j\leq k }}
 \Biggl[  \prod\sb{\substack{i\sb{j} \, \text{ s.t. }\\
i\sb{j}-i\sb{j-1}(>1) \\  2\leq j \leq k}} (m-2)(m-3)\sp{i\sb{j}-i\sb{j-1}-2} \Biggr] \cdot  (m-2)(m-3)\sp{n-i\sb{k}-1}  \\
& \quad + 2\sum\sb{\substack{\{  1 < i\sb{1} < i\sb{2} < i\sb{3} < \cdots < i\sb{k} < n  \} \\
\text{ s.t. } i\sb{j}-i\sb{j-1}>1, \;\;  2\leq j\leq k }}
 (m-2)(m-3)\sp{i\sb{1}+n-i\sb{k}-2}\cdot \Biggl[  \prod\sb{\substack{i\sb{j} \, \text{ s.t. }\\
i\sb{j}-i\sb{j-1}(>1) \\  2\leq j \leq k}} (m-2)(m-3)\sp{i\sb{j}-i\sb{j-1}-2} \Biggr] ,
\end{align*}
whereas
\[
v\sb{n, m}\sp{k, 1}=0 \qquad \text{ for } \lfloor n/2 \rfloor < k \leq n .
\]
\end{proposition}
\begin{proof}
We remark that
\[
\Hom(W\sb{n}, K\sb{m}) = \chi_{W\sb{n}}(m) = m(m-2)\biggl( (m-2)\sp{n-1}+(-1)\sp{n}\biggr)
\]
\cite[eq.\ (1.23), p.\ 13]{dkt}, and omit the rest of the proof for it is analogous to the one for Proposition \ref{prop:bwnkm1}.

\end{proof}

\begin{corollary}\label{cor:woddk31} For integer $n'\geq 1$,
\[
\hom(BW\sb{2n'+1}, K\sb{3}\sp{1}) = 0 .
\]
\end{corollary}
\begin{proof} Omitted.
\end{proof}

\section{Final remarks}\label{finrems}

We look forward to calculating numbers of homomorphisms from other graphs to the quasi-complete graphs.  We would also like to calculate partial profiles of other graphs. Here our plan will be to have for target graphs, graphs obtained from complete graphs by removing increasingly more edges from it. Furthermore, we would like to present our expressions in \emph{closed form} instead of the way the expressions for $\hom (P\sb{n}, K\sb{m}\sp{1})$, $\hom (C\sb{n}, K\sb{m}\sp{1})$, $\hom (BW\sb{n}, K\sb{m}\sp{1})$ and  $\hom (W\sb{n}, K\sb{m}\sp{1})$ are presented. The enumeration of homomorphisms of weighted graphs, which may be of significance in Statistical Mechanics of Exactly Solved Models \cite{fls}, is also in our plans.

Moreover, in this work, we also calculated the ``quasi''-chromatic ``polynomials'' of several graphs. As a matter of fact, we can regard our work here from the point of view of calculating $\hom(G, K\sb{m}\sp{1})$, for fixed $G$, which we call ``quasi''-chromatic ``polynomials''. This is also a potentially interesting line of work to pursue.

We plan to address these issues in future work.

\section{Acknowledgements}\label{ack}

The author acknowledges support by {\em Programa Operacional
``Ci\^{e}ncia, Tecnologia, Inova\c{c}\~{a}o''} (POCTI) of the {\em
Funda\c{c}\~{a}o para a Ci\^{e}ncia e a Tecnologia} (FCT)
cofinanced by the European Community fund FEDER. He thanks
the staff at IMPA and especially his host, Marcelo Viana, for
hospitality during his stay at this Institution. Finally, he thanks Louis H. Kauffman for suggestions.

\bigskip
\hrule
\bigskip

\noindent 2010 {\it Mathematics Subject Classification}:
Primary 05A15; Secondary 05C15.

\noindent \emph{Keywords: }
enumeration, sequence of integers, graph, graph homomorphism, complete graph,
quasi-complete graph, path, cycle, wheel, broken wheel, Lov\'asz vector of graph (profile of graph).

\bigskip
\hrule
\bigskip





\end{document}